\documentclass[preprint]{elsarticle}

\usepackage{amsmath}
\usepackage{amsthm}
\usepackage{amssymb}
\usepackage{bm}

\usepackage{tikz}
\usepackage{pgfplots}

\newtheorem{thm}{Theorem}

\newdefinition{rmk}{Remark}
\newproof{pf}{Proof}

\begin{document}

\begin{frontmatter}

\title{Factorized schemes of second-order accuracy for numerical solving unsteady problems\tnoteref{label1}}
\tnotetext[label1]{This work was supported by the Russian Foundation for Basic Research (project 14-01-00785)}

\author{P.N. Vabishchevich\corref{cor1}\fnref{lab1,lab2}}
\ead{vabishchevich@gmail.com}
\cortext[cor1]{Correspondibg author.}

\address[lab1]{Nuclear Safety Institute, Russian Academy of Sciences,
              52, B. Tulskaya, 115191 Moscow, Russia}

\address[lab2]{North-Eastern Federal University,
	      58, Belinskogo, 677000 Yakutsk, Russia}

\begin{abstract}

Schemes with the second-order approximation in time are considered for numerical solving the Cauchy problem 
for an evolutionary equation of first order with a self-adjoint operator. 
The implicit two-level scheme based on the Pad\'{e} polynomial approximation
is unconditionally stable. It demonstrates good asymptotic properties in time
and provides an adequate evolution in time for individual harmonics of the solution (has spectral mimetic stability).
In fact, the only drawback of this scheme is the necessity to solve 
an equation with an operator polynomial of second degree at each time level.
We consider modifications of these schemes, which are based on solving equations with operator polynomials of first degree.
Such computational implementations occur, for example,
if we apply the fully implicit two-level scheme (the backward Euler scheme).
A three-level modification of the SM-stable scheme is proposed. Its unconditional stability
is established in the corresponding norms. The emphasis is on the scheme, where the numerical algorithm
involves two stages, namely, the backward Euler scheme of first order at the first (prediction) stage 
and the following correction of the approximate solution using a factorized operator.
The SM-stability is established for the proposed scheme. To illustrate the theoretical results of the work, 
a model problem is solved numerically.

\end{abstract}

\begin{keyword}
Evolutionary equation of first order \sep the Cauchy problem \sep finite difference schemes \sep
Pad\'{e}  approximation \sep  SM-stability \sep factorized scheme

\PACS 02.30.Jr \sep 02.60.Lj \sep 02.70.Bf

\MSC 65J08 \sep  65M06  \sep 65M12

\end{keyword}

\end{frontmatter}

\section{Introduction} 
\label{sec:1}

In the practice of scientific computations for numerical solving transient boundary value problems,
focus is on computational algorithms of higher accuracy (see, e.g., \cite{HundsdorferVerwer2003,Gustafsson2008}).
Along with increasing the accuracy of approximation in space,
researchers try to improve the accuracy of approximation in time using primarily 
numerical methods developed for ordinary differential equations \cite{Ascher,LeVeque}. 
Concerning unsteady problems governed by partial differential equations,
we are interested in numerical methods to solve the Cauchy problem for stiff
systems of ordinary differential equations \cite{RakitskiiUstinovChernoruckii1979,HairerWanner2010,Butcher2008}.

Increasing of the accuracy of the approximate solution to time-dependent problems is achieved in different ways.
For two-level schemes, where the solution at two consecutive time levels is involved, 
polynomial approximations are  explicitly or implicitly employed for operators of difference schemes.
The Runge-Kutta methods \cite{Butcher2008,DekkerVerwer1984} are well-known examples of such schemes
widely used in modern computing practice. The main feature of multilevel schemes (multistep methods) results in
the approximation of time derivatives with a higher accuracy on a multi-point stencil. 
Multistep methods based on numerical backward differentiation formulas \cite{Gear1971}
should be mentioned as a typical example.

In the transition from continuous problems to discrete ones, we seek to inherit the most 
important properties of the differential boundary value problem under the consideration. 
For time-dependent problems, emphasis is on the main criterion to well-posedness of the grid problem, i.e.,
the stability of the difference solution with respect to small perturbations in
the initial data, boundary conditions, right-hand sides and coefficients of
equations \cite{Samarskii1989}.

Various classes of stable difference schemes can be used. Among stable difference schemes, we
search for a scheme that is optimal in sense of some additional criteria.
In the theory of difference schemes, the class of asymptotically stable
difference schemes \cite{SamarskiiGulin1973,SamarskiiVabishchevich1995} was highlighted
as the class of methods ensuring the correct behavior of the approximate solution at long time.
In the theory of numerical methods for solving systems of ordinary differential
equations \cite{Butcher2008,Gear1971}, the concept of $L$-stable methods was introduced, where
an asymptotic behavior of the approximate solution is also investigated at long time.

Spectral Mimetic (SM) properties of operator-difference schemes for numerical solving
the Cauchy problem for evolutionary equations of first order are  associated with the time-evolution of
individual harmonics of the solution. 
A study of spectral characteristics allows to select more acceptable approximations in time.
The concept of SM-stable difference schemes has been introduced in \cite{Vabischevich2010b}
for two-level schemes of higher accuracy for solving the Cauchy problem for evolutionary
first-order equations. The long-time behavior of individual harmonic of the approximate solution was studied
for problems with a self-adjoint operator using various implicit schemes, which are based on
the Pad\'{e} approximations \cite{BakerGraves-Morris1996} with the corresponding operator exponential.
In \cite{Vabischevich2010b,Vabischevich2011a}, it was shown that the best way to obtain SM-properties for 
such problems is to apply polynomial approximations.

Among two-level implicit schemes of higher accuracy order for evolutionary  equations with skew-symmetric operators
that were obtained on the basis of the Pad\'{e} approximations, symmetric schemes \cite{Vabischevich2011} 
demonstrate some advantages. For evolutionary equations with general operators, it is possible 
to construct additive schemes (splitting schemes) via splitting the problem operator into
self-adjoint and skew-symmetric components \cite{Vabischevich2012}.

The analysis of schemes for problems with self-adjoint operators conducted in the mentioned above papers
allowed to highlight the class of SM-stable difference schemes based on the Pad\'{e} polynomial approximation.
The implicit two-level scheme of second order approximation in time constructed using 
the Pad\'{e}  polynomial approximation is unconditionally SM-stable.
The computational implementation of this scheme involves the solution of  an
equation with an operator polynomial of second degree at each time step.
It seems reasonable to modify these schemes in order to employ more simple equations with an operator polynomial of first degree. For example, we can use the fully implicit (backward Euler) two-level scheme.
In this case, at each time level, we solve the problem with a factorized operator,
where each factor is an operator polynomial of first degree.
When considering systems of ordinary differential equations, in this case, we speak of
diagonally implicit Runge-Kutta (DIRK methods) \cite{Butcher2008,DekkerVerwer1984}.
Factorized SM-stable schemes were investigated in the work \cite{Vabischevich2010c}.
It was found that among the schemes with maximal accuracy, in fact, there is no SM-stable scheme.

In the present paper, we construct factorized schemes of second-order approximation in time 
to the Cauchy problem for the evolutionary equation of first order with a self-adjoint operator. 
They are based on a modification of the implicit two-level scheme based on the Pad\'{e}  polynomial approximation, 
which has the second-order approximation in time, is unconditionally stable, 
demonstrates good asymptotic properties in time and adequate long-time behavior
of individual harmonics of the solution (spectral mimetic stability).

In the first modification, instead of two-level scheme, we employ three-level scheme.
Here a part of the approximation for the time derivative is taken from the previous time level.
A similar idea we already used  (see, e.g., \cite{Vabischevich2015})
to increase the accuracy of explicit-implicit schemes. The unconditional stability of the three-level
scheme is shown in the corresponding norms. A more interesting modification
employs the two-stage calculations, namely, the backward Euler scheme of first order at the prediction stage
and the following correction of the approximate solution via solving the problem with a factorized operator.
This scheme is unconditionally SM-stable.

The paper is organized as follows. 
In Section 2, we consider unconditionally SM-stable scheme of second-order approximation in time 
for solving the Cauchy problem for the evolutionary equation of first order. 
The factorized scheme based on  the three-level modification is investigated in Section 3. 
Section 4 is the core of our work. Here we study the SM-stable factorized scheme. 
Numerical experiments for a model problem with a one-dimensional parabolic 
equation are presented in Section 5.

\section{SM-stable scheme of second order} 

Let us consider a finite-dimensional real Hilbert space $H$, where the scalar product and the norm 
are $(\cdot,\cdot)$ and $\|\cdot\|$, respectively. 
Suppose that  $u(t)$ ($0 \leq t \leq T$, $T > 0$) is defined as the solution of the Cauchy problem 
for the homogeneous evolutionary first-order equation:
\begin{equation}\label{1}
  \frac{d u}{d t} + A \, u = 0,
  \quad 0 < t \leq T ,
\end{equation}
\begin{equation}\label{2}
  u(0) = u^0 . 
\end{equation}
In equation (\ref{1}), $A$ is linear and, for simplicity, independent of $t$ (constant) 
operator acting from $H$ into $H$ ($ A \, d/dt = d/dt \, A$). In addition, assume that the operator is  
self-adjoint positive definite: 
\begin{equation}\label{3}
  A = A^* \geq \delta E,
  \quad \delta \geq 0,
\end{equation}
where $E$ is the identity operator.

For problem (\ref{1})--(\ref{3}), the following elementary estimate 
for stability of the solution with respect to the initial data holds:
\begin{equation}\label{4}
  \|u(t)\| \leq 
  \exp(-\delta t) \|u^0\|.
\end{equation}
SM-properties for approximations in time are connected \cite{Vabischevich2010b} 
with inheriting the basic properties of the differential problem.

To solve numerically the problem (\ref{1}), (\ref{2}), 
we employ two-level difference schemes. Define the uniform grid in time with a step $\tau$ as
\[
  \overline{\omega}_\tau =
  \omega_\tau\cup \{T\} =
  \{t^n=n\tau,
  \quad n=0,1,...,N,
  \quad \tau N=T\} 
\]
and denote $y^n = y(t^n)$, where $t^n = n \tau$. We want to pass from a time level $t^n$ 
to the next time level $t^{n+1}$. The exact solution has the form
\begin{equation}\label{5}
  u(x,t^{n+1}) =  \exp(- \tau A) u(x,t^{n})  .
\end{equation}

The two-level scheme corresponding to (\ref{5})
may be written in the form
\begin{equation}\label{6}
y^{n+1} = S y^{n}, \qquad  n = 0,1, ... ,
\end{equation}
where $S$ is the transition operator from the current time level to the next one. 
In the general case, the operator $S$ may depend on~$n$.
We restrict the consideration to the difference approximations with respect to time 
for the homogeneous problem (\ref{1}), (\ref{2}) that lead to the transition operator
\begin{equation}\label{7}
 S = s(\tau A ),
\end{equation}
where $s(z)$ is called a stability function. 

The use of additional information about properties of a problem operator 
allows to specify requirements for the choice of approximations in time.
Taking into account the inequality (\ref{3}),
we get a more accurate stability estimate 
\[
  \| u(x,t^{n+1}) \| \leq \exp(-\delta \tau) \| u(x,t^{n}) \|.
\]
The inheritance of this property of the problem
(\ref{1}), (\ref{2}) means that the following estimate must be true:
\begin{equation}\label{8}
  \| y^{n+1} \| \leq \rho \| y^{n} \|,
  \quad \|s(\tau A ) \| \leq \rho ,
  \quad \rho  =  \exp(-\delta \tau ) .  
\end{equation}
In this case, we say \cite{SamarskiiGulin1973} that this operator-difference (semi-discrete) scheme is $\rho$-stable.

The properties discussed above concern the bounds of the spectrum of the transition operator. We 
want also to study some other qualitative indicators of the approximate solution  behavior.
In solving unsteady problems, we focuses on the long-time asymptotic behavior of the solution.
For the model problem under examination, the solution decreases to zero as $t\rightarrow \infty$. 
If the approximate solution preserves this property, we say that the approximate solution 
is asymptotically stable (with respect to time). In the case of the Cauchy problem for systems of ODEs, 
this property of the approximate solution is called $L$-stability. If
\begin{equation}\label{9}
  \lim_{\theta \rightarrow \infty } s(\theta ) = 0 ,
\end{equation}
then the stable ($\rho$-stable) difference scheme (\ref{6}) is said to be asymptotically stable.

For the linear problem (\ref{1}), (\ref{2}) with a self-adjoint positive
definite operator $A$, the solution may be written as the superposition of individual harmonics, which 
are associated with their own eigenvalues. In \cite{Vabischevich2010b},
the choice of approximations in time is subjected to the requirement of
the appropriate behavior in time of all points of the spectrum.
In this case, we introduce SM-property (Spectral Mimetic) of schemes, i.e., such a scheme is said to be SM-stable.
A difference scheme is called SM-stable if it is $\rho$-stable and asymptotically stable.
The additional requirement is the spectral monotonicity, namely,
the function $s(\theta)$ is monotonically decreasing. 
This means that the harmonics with higher indexes decay more rapidly than the harmonics with lower indexes.

Two-level schemes of higher-order approximations for time-dependent linear problems can be conveniently 
constructed on the basis of the Pad\'{e} approximations of the operator (matrix) exponential function $\exp(- \tau A)$. 
Other approaches are discussed in \cite{Higham2008,MolerLoan2003}. In the case of nonlinear systems of
ODEs, such approximations correspond to various variants of Runge-Kutta methods 
\cite{HundsdorferVerwer2003,HairerWanner2010,Butcher2008}.

The Pad\'{e} approximation of the function $\exp(-z)$ is
\begin{equation}\label{10}
  \exp(-z) = R_{lm}(z) + \mathcal{O} (z^{l+m+1}),
  \quad R_{lm}(z) \equiv \frac{P_{lm}(z)}{Q_{lm}(z)},
\end{equation}
where $P_{lm}(z)$ and $Q_{lm}(z)$ are polynomials of degrees $l$ and $m$, respectively. 
These polynomials have \cite{BakerGraves-Morris1996} the form
\[
  P_{lm}(z) = \frac{l!}{(l+m)!} \sum_{k=0}^{l} 
  \frac{(l+m-k)!}{k! (l-k)!} (-z)^k,
\]
\[
  Q_{lm}(z) = \frac{m!}{(l+m)!} \sum_{k=0}^{m} 
  \frac{(l+m-k)!}{k! (m-k)!} z^k .
\]

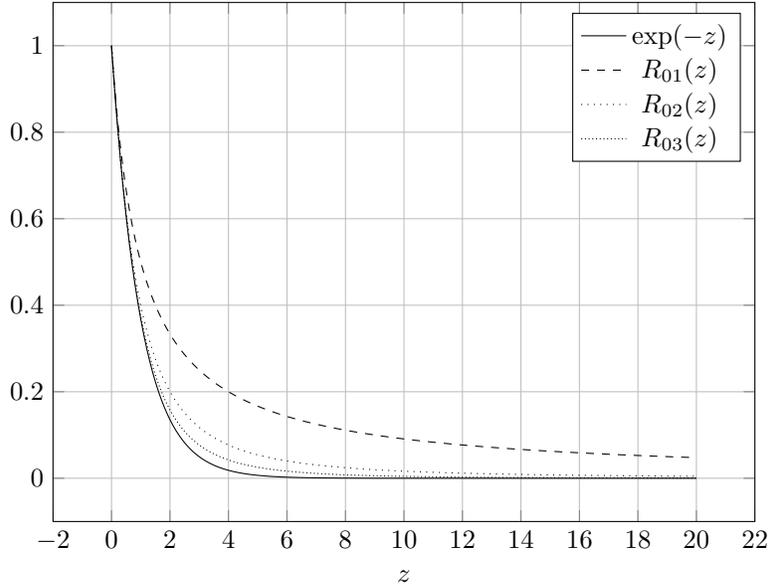
\begin{figure}[ht] 
  \begin{center}
    \begin{tikzpicture}
    \begin{axis}[height=0.7\textwidth, width=0.9\textwidth, grid=major,
    domain=0:20, samples=200,
    xlabel=$z$]
    \addplot[color=black,style=solid] {exp(-x)};
    \addlegendentry{$\exp(-z)$}
    \addplot[color=black,style=dashed] {1/(1+x)};
    \addlegendentry{$R_{01}(z)$}
    \addplot[color=black,style=dotted]{2/(2+2*x+x^2)};
    \addlegendentry{$R_{02}(z)$}
    \addplot [color=black,style=densely dotted]{6/(6+6*x+3*x^2+x^3)};
    \addlegendentry{$R_{03}(z)$}
    \end{axis}
    \end{tikzpicture}
    \caption{SM-stable Pad\'{e} approximations} 
	\label{f-1}
  \end{center}
\end{figure}

Only the schemes with $l = 0$, i.e.,
\[
  R_{0m}(z) \equiv \frac{1}{Q_{0m}(z)},
  \quad P_{0m}(z) = 1
\]
are SM-stable. In this case, the two-level difference scheme for the problem (\ref{1}), (\ref{2}) is
\begin{equation}\label{11}
  Q_{0m}(\tau \Lambda) \frac{y^{n+1} -y^n}{\tau} + 
  \frac{1}{\tau }(Q_{0m}(\tau \Lambda) - E) y^{n} = 0, 
  \quad  n = 0,1, ... ,  
\end{equation}
where the function
\[
  Q_{0m}(z) = \sum_{k=0}^{m} \frac{1}{k!} z^k 
\]
is a truncated Taylor series for $\exp(z)$. Figure~\ref{f-1} presents the corresponding stability 
functions for $m=1,2,3$.

When considering the Cauchy problem with an inhomogeneous right-hand side, 
instead of (\ref{1}), we solve the equation 
\begin{equation}\label{12}
  \frac{d u}{d t} + A \, u = f(t),
  \quad 0 < t \leq T ,
\end{equation}
The SM-stable scheme of first order ($m=1$ in (\ref{11})) may be written in the form
\begin{equation}\label{13}
 (E + \tau A)\frac{y^{n+1}-y^{n}}{\tau} + A y^{n} = \varphi^n .
\end{equation} 
Taking into account that the exact solution of (\ref{2}), (\ref{12}) satisfies 
the representation
\[
\begin{split}
  u(t^{n+1}) & = \exp\left (- \tau  \sum_{\alpha = 1}^{p} A_\alpha \right )  u(t^{n}) \\
  & +
  \int_{t^n}^{t^{n+1}} \exp\left (-(t^{n+1}-\theta)  
  \sum_{\alpha = 1}^{p} A_\alpha \right ) f(\theta) d \theta ,  
\end{split}
\]
for the right-hand side of equations (\ref{13}), we can take $\varphi^n = f(t^{n+1/2})$.

The object of our study is the factorized scheme of second order ($m=2$ in (\ref{11})), where
\begin{equation}\label{14}
 \left (E + \tau A + \frac{1}{2} \tau^2 A^2 \right )\frac{y^{n+1}-y^{n}}{\tau} + 
 \left (A + \frac{1}{2} \tau A^2 \right ) y^{n} = \varphi^n .
\end{equation}
For the right-hand side, it seems reasonable to use the approximation
\[
 \varphi^n = \left (E + \frac{1}{2} \tau A \right )f(t^{n+1/2}) .
\]
The main drawback of the scheme (\ref{14}) is associated with solving the problem
\[
 \left (E + \tau A + \frac{1}{2} \tau^2 A^2 \right ) y^{n+1} = \chi^n
\]
at each time level, where we have $A^2$.
The computational implementation of the backward Euler scheme (\ref{13}) involves 
the solution of the essentially simpler problem $(E + \tau A) y^{n+1} = \chi^n$.
A fundamental simplification in solving the problem at a new time level is connected
with using the factorized scheme, where we solve the problem
\begin{equation}\label{15}
 (E + \sigma_1\tau A) (E + \sigma_2 \tau A) y^{n+1} = \chi^n 
\end{equation}
with positive constants $\sigma_\alpha, \ \alpha = 1,2$. 

\section{Three-level factorized scheme} 

To construct factorized schemes,  we restrict our study to the case with the identical weighting parameters 
in (\ref{15}), i.e., the solution at the new time level is evaluated from the equation
\[
 (E + \sigma\tau A)^2 y^{n+1} = \chi^n .
\]
Let us rewrite the basic scheme (\ref{14}) in the form
\begin{equation}\label{16}
\begin{split}
 (E + \sigma\tau A)^2 \frac{y^{n+1}-y^{n}}{\tau} & +
 \left (E + \tau A + \frac{1}{2} \tau^2 A^2 - (E + \sigma\tau A)^2 \right )\frac{y^{n+1}-y^{n}}{\tau} \\
 & +  \left (A + \frac{1}{2} \tau A^2 \right ) y^{n} = \varphi^n .
\end{split}
\end{equation}
Factorized schemes are based on some modification of the second term in the left-hand side of (\ref{16}). 
In view of
\[
 E + \tau A + \frac{1}{2} \tau^2 A^2 - (E + \sigma\tau A)^2 =
 (1-2 \sigma) \tau A + \left (\frac{1}{2} - \sigma^2 \right ) \tau^2 A^2 ,
\]
to preserve the second-order approximation, it is sufficient to approximate the difference derivative 
$(y^{n+1}-y^{n}) / \tau$ with the first order with respect to $\tau$.

Consider the case of modification, where
\[
 \frac{y^{n+1}-y^{n}}{\tau} \longrightarrow \frac{y^{n}-y^{n-1}}{\tau} .
\]
Such a transition from a two-level scheme to a three-level scheme was used in 
\cite{Vabischevich2015} in the construction of explicit schemes for parabolic and hyperbolic equations.

Instead of (\ref{14}), we apply the scheme
\begin{equation}\label{17}
\begin{split}
 (E + \sigma\tau A)^2 \frac{y^{n+1}-y^{n}}{\tau} & +
 \left (E + \tau A + \frac{1}{2} \tau^2 A^2 - (E + \sigma\tau A)^2 \right )\frac{y^{n}-y^{n-1}}{\tau} \\
 & +  \left (A + \frac{1}{2} \tau A^2 \right ) y^{n} = \varphi^n .
\end{split}
\end{equation} 
Taking into account that
\[
 \frac{y^{n+1}-y^{n}}{\tau} = \frac{y^{n+1}-y^{n-1}}{2\tau} + \frac{\tau }{2} \frac{y^{n+1} - 2y^{n}+y^{n-1}}{\tau^2} ,
\] 
\[
 \frac{y^{n}-y^{n-1}}{\tau} = \frac{y^{n+1}-y^{n-1}}{2\tau} - \frac{\tau }{2} \frac{y^{n+1} - 2y^{n}+y^{n-1}}{\tau^2} ,
\]
rewrite it in the canonical form \cite{Samarskii1989} for three-level operator-difference schemes:
\begin{equation}\label{18}
 B \frac{y^{n+1}-y^{n}}{\tau}  + D \frac{y^{n+1} - 2y^{n}+y^{n-1}}{\tau^2} + \widetilde{A} y^{n} = \varphi^n .
\end{equation} 
For the operators $\widetilde{A}, B, D$, we have
\[
 \widetilde{A} = A + \frac{1}{2} \tau A^2,
\] 
\[
 B = E + \tau A + \frac{1}{2} \tau^2 A^2 ,
\] 
\[
 D = \frac{\tau }{2} \left ( E + (4 \sigma - 1) \tau A + \left (2 \sigma^2 -  \frac{1}{2} \right )  \tau^2 A^2 \right ) .
\] 

The necessary and sufficient conditions for the stability of the three-level scheme (\ref{18}) 
with constant self-adjoint  operators $\widetilde{A}, B, D$, where $B \geq 0$ and $\widetilde{A} > 0$,
are formulated (see \cite{Samarskii1989,SamarskiiGulin1973,SamarskiiMatusVabishchevich}) in the form of the inequality
\begin{equation}\label{19}
 D > \frac{\tau^2}{4} \widetilde{A} . 
\end{equation} 
We have
\[
 G = D - \frac{\tau^2}{4} \widetilde{A} =
 \frac{\tau }{2} \left ( E + \left (4 \sigma - \frac{3}{2} \right ) \tau A + \left (2 \sigma^2 -  \frac{3}{4} \right )  \tau^2 A^2 \right ) .
\] 
The inequality (\ref{19}) holds for $\sigma \geq \sqrt{3/8}$.

\begin{thm} 
\label{t-1} 
The three-level factorized scheme (\ref{17}) 
is unconditionally stable under the restriction $\sigma \geq \sqrt{3/8}$. The following a priori estimate holds:
\begin{equation}\label{20}
  \mathcal{E}_{n+1} \leq 
  \mathcal{E}_{n} +
  \frac{\tau}{2}  \|\varphi^n \|^2_{B^{-1}} ,
\end{equation}
where
\[
 \mathcal{E}_{n+1} =
 \left \| \frac{y^{n+1} + y^{n}}{2 } \right \|^2_{\widetilde{A}}+
 \left \| \frac{y^{n+1} - y^{n}}{\tau } \right \|^2_{G} .
\] 
\end{thm} 
 
\begin{pf}
By
\[
  y^{n} = 
  \frac{1}{4} (y^{n+1} + 2 y^{n} + y^{n-1}) -
  \frac{1}{4} (y^{n+1} - 2 y^{n} + y^{n-1}) ,
\] 
we can rewrite (\ref{18}) as 
\begin{equation}\label{21}
  B \frac{y^{n+1} - y^{n-1}}{2\tau }  +
  G \frac{y^{n+1} - 2 y^{n} + y^{n-1}}{\tau^2} +
  \widetilde{A} \frac{y^{n+1} - 2 y^{n} + y^{n-1}}{4} =
  \varphi^{n} .
\end{equation} 
Let  
\[
  v^{n} = \frac{1}{2} (y^{n} + y^{n-1}),
  \quad w^{n} = \frac{y^{n} - y^{n-1}}{\tau} ,
\]
then (\ref{21}) can be written in the form 
\begin{equation}\label{22}
  B \frac{w^{n+1} + w^{n}}{2} 
  + G \frac{w^{n+1} - w^{n}}{\tau } +
  \frac{1 }{2} \widetilde{A}  ( v^{n+1} + v^{n})  =
  \varphi^{n} .
\end{equation}
Multiplying scalarly both sides of (\ref{22}) by  
\[
  2 (v^{n+1} - v^{n}) =
  \tau (w^{n+1} + w^{n}) ,
\]
we get the equality
\begin{equation}\label{23}
\begin{split}
  \frac{\tau}{2} 
  ( B (w^{n+1} & + w^{n}), w^{n+1} + w^{n}) +
  ( G (w^{n+1} - w^{n}), w^{n+1} + w^{n}) 
  \\
  & +
  ( \widetilde{A} (v^{n+1} + v^{n}), v^{n+1} - v^{n}) =
  \tau (\varphi^{n}, w^{n+1} + w^{n} ) .
\end{split}
\end{equation}
To estimate the right-hand side, we use the inequality 
\[
  ( \varphi^{n}, w^{n+1} + w^{n} ) \leq 
  \frac{1 }{2} 
  ( B (w^{n+1} + w^{n}), w^{n+1} + w^{n}) +
  \frac{1}{2} 
  (B^{-1} \varphi^n, \varphi^n) .
\]
This makes it possible to get from (\ref{23}) the inequality 
\begin{equation}\label{24}
  \mathcal{E}_{n+1} \leq 
  \mathcal{E}_{n} +
  \frac{\tau}{2} 
  (B^{-1} \varphi^n, \varphi^n) ,
\end{equation}
where we use the notation 
\[
  \mathcal{E}_{n} = 
  ( \widetilde{A} v^{n}, v^{n})
  + ( G w^{n}, w^{n}) .
\]
The inequality (\ref{24}) is the desired a priori estimate (\ref{20}). 
\end{pf}

\section{SM-stable factorized schemes} 

The second possibility to modify the scheme (\ref{16}) is connected 
with the solution of an auxiliary problem. In this case
\[
 \frac{y^{n+1}-y^{n}}{\tau} \longrightarrow \frac{\widetilde{y}^{n+1}-y^{n}}{\tau} ,
\]
in calculating the auxiliary quantity $\widetilde{y}^{n+1}$, which approximates $y^{n+1}$.
This stage of predicting the solution at the new time level can be naturally implemented 
on the basis of SM-stable scheme of first order.  According to (\ref{13}), we put
\begin{equation}\label{25}
 (E + \tau A)\frac{\widetilde{y}^{n+1}-y^{n}}{\tau} + A y^{n} = \widetilde{\varphi}^n ,
\end{equation}
where, for example, $\widetilde{\varphi}^n = f(t^{n+1/2})$. At the correction stage,
we apply the factorized scheme
\begin{equation}\label{26}
\begin{split}
 (E + \sigma\tau A)^2 \frac{y^{n+1}-y^{n}}{\tau} & +
 \left (E + \tau A + \frac{1}{2} \tau^2 A^2 - (E + \sigma\tau A)^2 \right )\frac{\widetilde{y}^{n+1}-y^{n}}{\tau} \\
 & +  \left (A + \frac{1}{2} \tau A^2 \right ) y^{n} = \varphi^n .
\end{split}
\end{equation}
Let us formulate the conditions for SM-stability of the factorized scheme (\ref{25}), (\ref{26}).

\begin{thm} 
\label{t-2} 
The  factorized scheme (\ref{25}), (\ref{26})
is unconditionally SM-stable under the restriction $\sigma \geq 1/\sqrt{2}$. 
\end{thm} 
 
\begin{pf}
From (\ref{25}), we have
\[
\frac{\widetilde{y}^{n+1}-y^{n}}{\tau}  = - ((E + \tau A)^{-1} (A y^{n} - \widetilde{\varphi}^n ) . 
\] 
Substitution in (\ref{26}) results in
\begin{equation}\label{27}
\begin{split}
 (E + \tau A) & (E + \sigma\tau A)^2 \frac{y^{n+1}-y^{n}}{\tau} -
 \left ((1-2\sigma)\tau A + \left ( \frac{1}{2} - \sigma^2 \right ) \tau^2 A^2 \right ) A y^{n} \\
 & + (E + \tau A) \left (E + \frac{1}{2} \tau A \right ) A y^{n} = \widetilde{\psi}^n .
\end{split}
\end{equation} 
The two-level scheme (\ref{27}) may be written as
\[
 y^{n+1} = S y^{n} + \tau \psi^{n} .
\]
For the stability function $s(z)$, we have the representation
\begin{equation}\label{28}
 s(z) = \frac{1 + 2\sigma z + (\sigma^2 - 0.5 ) z^2}{1 + (1+2\sigma) z + (\sigma^2+2\sigma) z^2 + \sigma^2 z^3} .
\end{equation}
It is easy to see that the condition of asymptotic stability (\ref{9}) holds for (\ref{28}).
The most important condition of SM-stability that is associated with  monotonous decreasing the function 
$s(z)$ with $z \geq 0$ is satisfied for $\sigma \geq 1/\sqrt{2}$. 
This proves the theorem.
\end{pf}

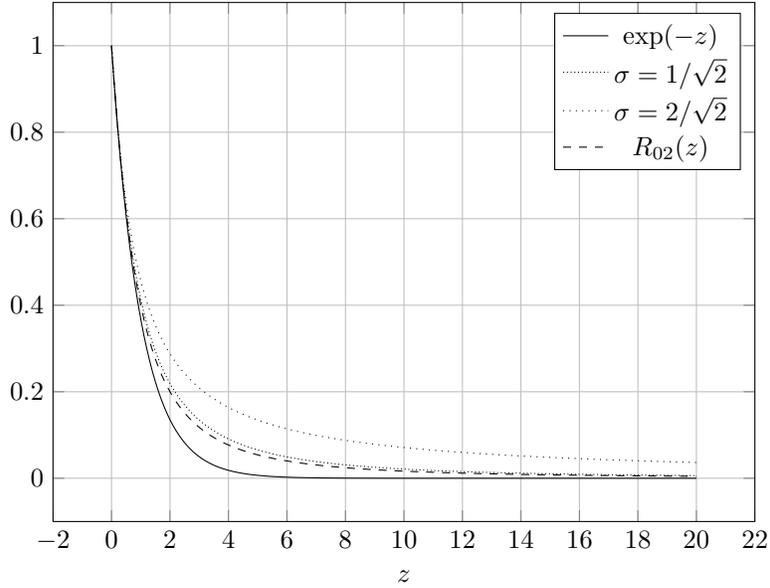
\begin{figure}[ht] 
  \begin{center}
    \begin{tikzpicture}
    \begin{axis}[height=0.7\textwidth, width=0.9\textwidth, grid=major,
    domain=0:20, samples=200,
    xlabel=$z$]
    \addplot[color=black,style=solid] {exp(-x)};
    \addlegendentry{$\exp(-z)$}
    \addplot[color=black,style=densely dotted] {(1+2*0.707*x)/(1+(1+2*0.707)*x + 0.707*(0.707+2)*x^2 + 0.707^2*x^3)};
    \addlegendentry{$\sigma = 1/\sqrt{2}$}
    \addplot[color=black,style=dotted] {(1+4*0.707*x + (4*0.707^2-0.5)*x^2)/(1+(1+4*0.707)*x + 2*0.707*(2*0.707+2)*x^2 + 4*0.707^2*x^3)};
    \addlegendentry{$\sigma = 2/\sqrt{2}$}
    \addplot[color=black,style=dashed]{2/(2+2*x+x^2)};
    \addlegendentry{$R_{02}(z)$}
    \end{axis}
    \end{tikzpicture}
    \caption{SM-stable factorized scheme} 
	\label{f-2}
  \end{center}
\end{figure}

The limiting value $\sigma = 1/\sqrt{2}$ provide the best approximation properties. 
It can be observed in Figure~\ref{f-2}. 
Note that the SM-stable scheme (\ref{25}), (\ref{26}) is based on the factorization 
into three operators of a simpler structure, whereas in the three-level scheme (\ref{17}), we have only two factors. 
The main advantage of the scheme (\ref{25}), (\ref{26}) is in its absolute SM-stability.

\section{Numerical experiments} 

Now we will illustrate the possibilities of the above factorized schemes of second-order accuracy 
in time on predictions of a model parabolic problem. Consider a boundary value problem for the equation
\begin{equation}\label{29}
 \frac{\partial u}{\partial t} - \frac{\partial^2 u}{\partial x^2} = f(x,t),
 \quad 0 < x < 1,
 \quad 0 < t \leq T . 
\end{equation} 
The equation (\ref{29}) is supplemented with the boundary and initial conditions:
\begin{equation}\label{30}
 u(0,t) = 0,
 \quad u(1,t) = 0, 
 \quad 0 < t \leq T ,
\end{equation} 
\begin{equation}\label{31}
 u(x,0) = 0,
 \quad 0 < x < 1 .
\end{equation} 
The problem (\ref{29})--(\ref{31}) is considered with $T = 0.5$, when the exact solution is
\[
 u(x,t) = x(1-x) (1-\exp(- 5 t)) . 
\] 

To solve numerically the problem (\ref{29})--(\ref{31}), we 
introduce a difference approximation in space. 
Within the unit interval, we introduce a uniform grid with step $h$:
\[
 \overline{\omega} = \omega \cup \partial \omega = \{ x \ | \ x = x_i = ih, \ i = 0,1, ..., M, \ M h = 1 \} ,
\]  
where $\omega$ is the set of internal nodes ($i = 1, 2, ..., M-1$) of the grid.
For $u(x,t), \ x \in \omega$, we have the equation
\[
 \frac{d y}{d t} + A y = f(x,t),
 \quad x \in \omega .
\] 
Here the discrete operator $A$  is defined as follows:
\[
 A v = \left \{
 \begin{array}{lll}
  {\displaystyle \frac{2 v(x) - v(x+h)}{h^2} } ,  &  x = x_1 , \\
  {\displaystyle \frac{2 v(x) - v(x-h) - v(x+h)}{h^2} } ,  &  x = x_i , \ i = 2, 3, ..., M-2 ,\\
  {\displaystyle \frac{2 v(x) - v(x-h)}{h^2} } ,  &  x = x_{M-1} . \\
\end{array}
\right . 
\]
The accuracy of the approximate solution is estimated in the grid norm of the space $L_2(\omega)$, 
where the scalar product and the norm are defined as
\[
 (v,w) = \sum_{x \in \omega} v(x) w(x) h,
 \quad \|v\| = (v, v)^{1/2} . 
\] 
Thus, for the error of the approximate solution $y^n - u(x,t^n), \ x \in \omega$, we study
$\varepsilon(t) = \|y^n - u(x,t)\|, \ t = t^n$.

In our numerical experiments, emphasis is on the accuracy in time of the difference schemes under consideration.
Because of this, we present predictions conducted using the same spatial grid $M = 10$. 
The calculations are performed using different grids in time, namely, $N = 10, 20, 40$.

\begin{figure}[htp]
  \begin{center}
    \includegraphics[scale = 0.5] {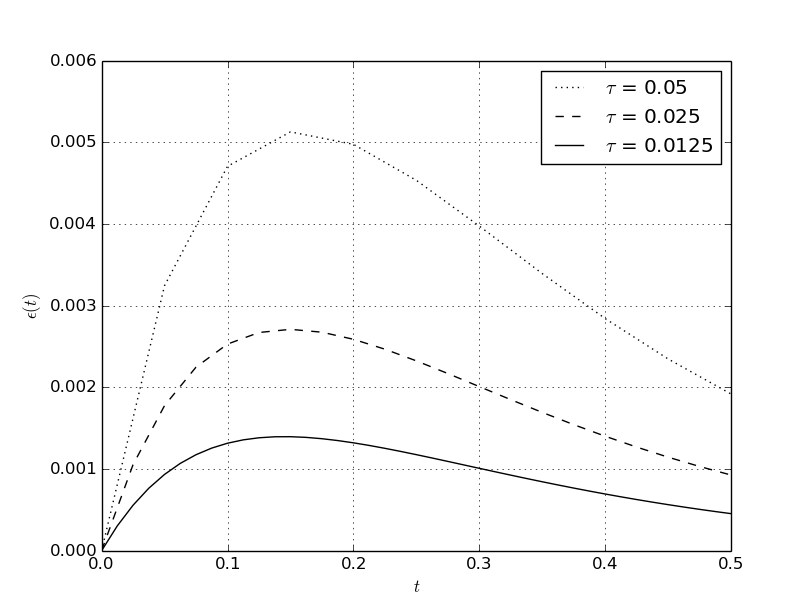}
	\caption{Backward Euler scheme}
	\label{f-3}
  \end{center}
\end{figure}

\begin{figure}[htp]
  \begin{center}
    \includegraphics[scale = 0.5] {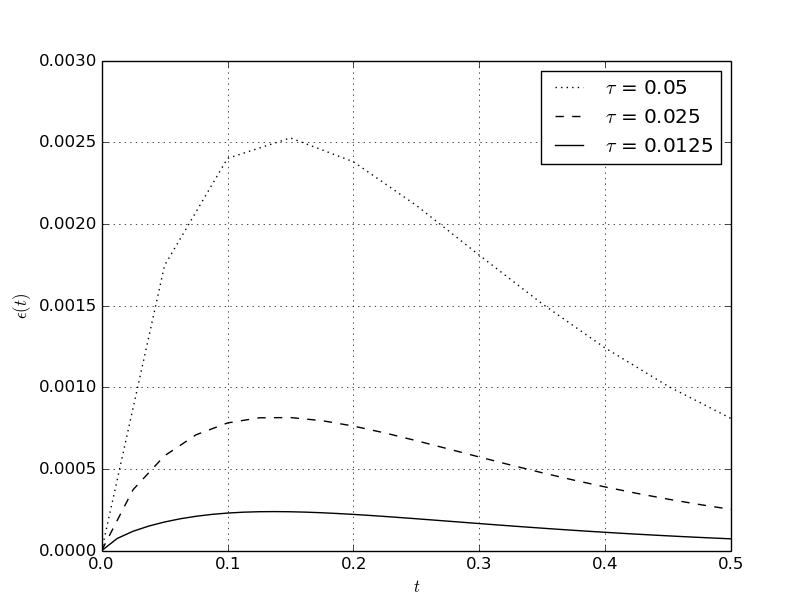}
	\caption{Crank-Nicolson scheme}
	\label{f-4}
  \end{center}
\end{figure}

The accuracy of conventional two-level difference schemes is presented in Figure~\ref{f-3} and Figure~\ref{f-4}. 
Figure~\ref{f-3} demonstrates the results for the  backward Euler scheme (\ref{13}) (the scheme of first-order accuracy). 
Similar numerical data for the Crank-Nicolson scheme
\[
  \left (E + \frac{1}{2} \tau A \right )\frac{y^{n+1}-y^{n}}{\tau} + A y^{n} = \varphi^n 
\] 
are depicted in Figure~\ref{f-4} (the scheme of second-order accuracy).  

\begin{figure}[htp]
  \begin{center}
    \includegraphics[scale = 0.5] {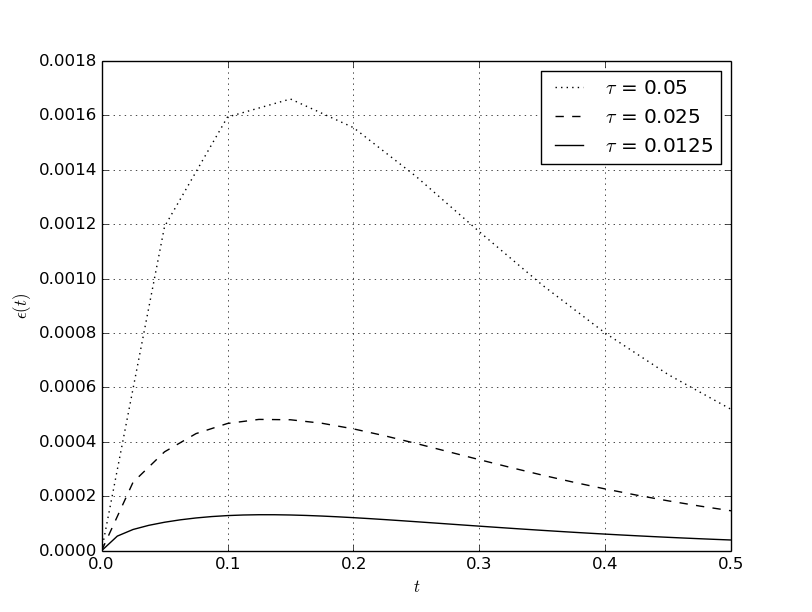}
	\caption{SM-stable scheme of second-order accuracy}
	\label{f-5}
  \end{center}
\end{figure} 

The main object of our study is the modified SM-stable scheme of second-order accuracy (\ref{14}).
The accuracy of this scheme is shown in Figure~\ref{f-5}.  
The approximate solution of the model problem (\ref{29})--(\ref{31}) obtained via 
the scheme (\ref{14}) demonstrates higher accuracy 
in comparison with the Crank-Nicolson scheme (compare Figure~\ref{f-5} and Figure~\ref{f-4}).

\begin{figure}[htp]
  \begin{center}
    \includegraphics[scale = 0.5] {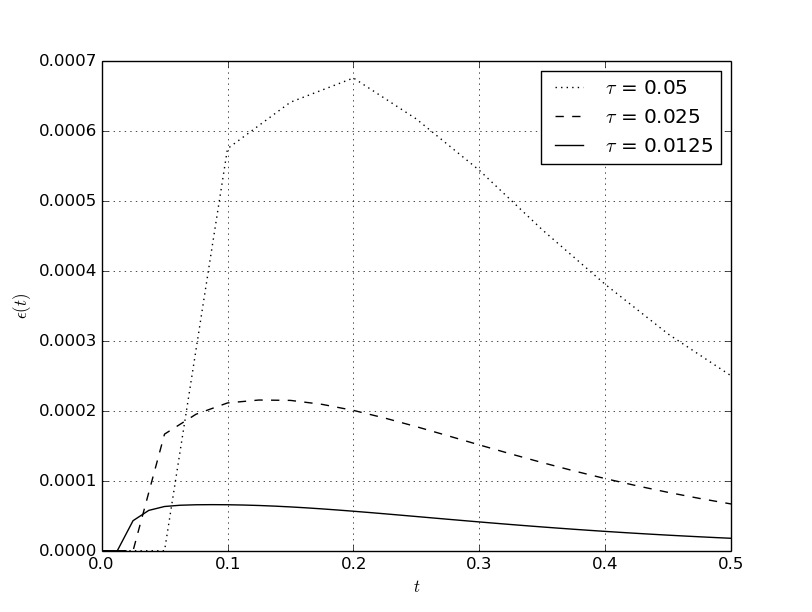}
	\caption{Three-level factorized scheme: exact $y^1$}
	\label{f-6}
  \end{center}
\end{figure}

\begin{figure}[htp]
  \begin{center}
    \includegraphics[scale = 0.5] {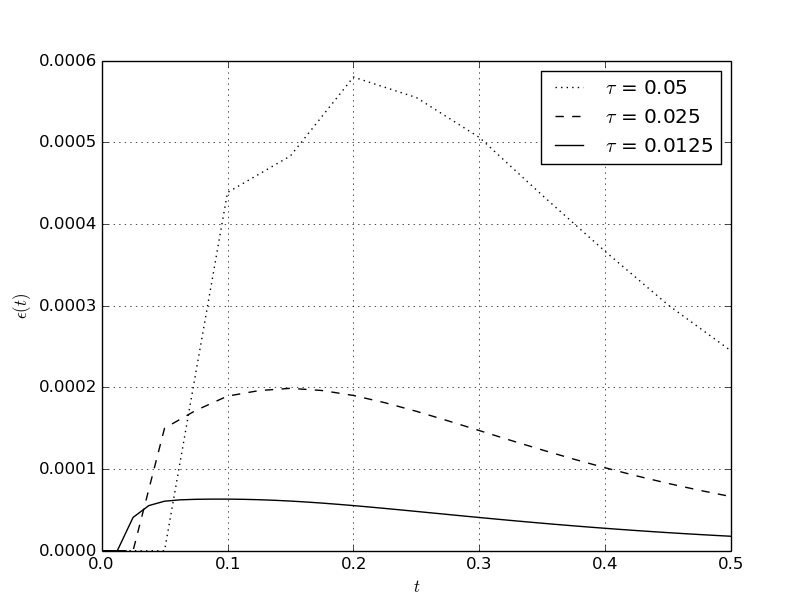}
	\caption{Three-level factorized scheme: the Crank-Nicolson scheme for evaluating $y^1$}
	\label{f-7}
  \end{center}
\end{figure}

When using the three-level factorized scheme (\ref{17}),   
we must preliminarily calculate the approximate solution at the first time level $y^1$.
Figure~\ref{f-6} demonstrates accuracy of predictions, where the exact $y^1$ is used. 
In practical calculations, we should focus on the schemes of the second-order accuracy for 
evaluating $y^1$. Thus, it seems natural to use the Crank-Nicolson scheme.  
In this case, the accuracy of calculations for the model problem (\ref{29})--(\ref{31}) is shown 
in Figure~\ref{f-7}.


\begin{figure}[htp]
  \begin{center}
    \includegraphics[scale = 0.5] {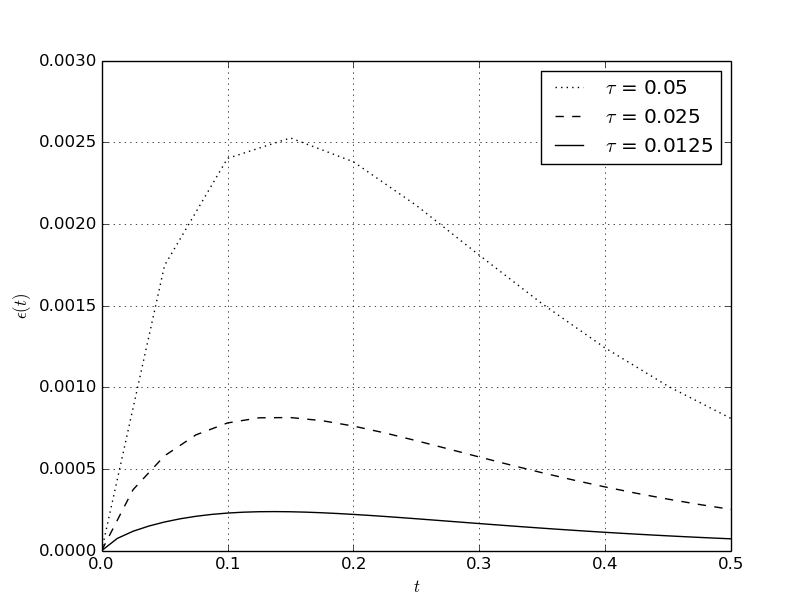}
	\caption{SM-stable factorized scheme: $\sigma = 1/\sqrt{2}$}
	\label{f-8}
  \end{center}
\end{figure}

\begin{figure}[htp]
  \begin{center}
    \includegraphics[scale = 0.5] {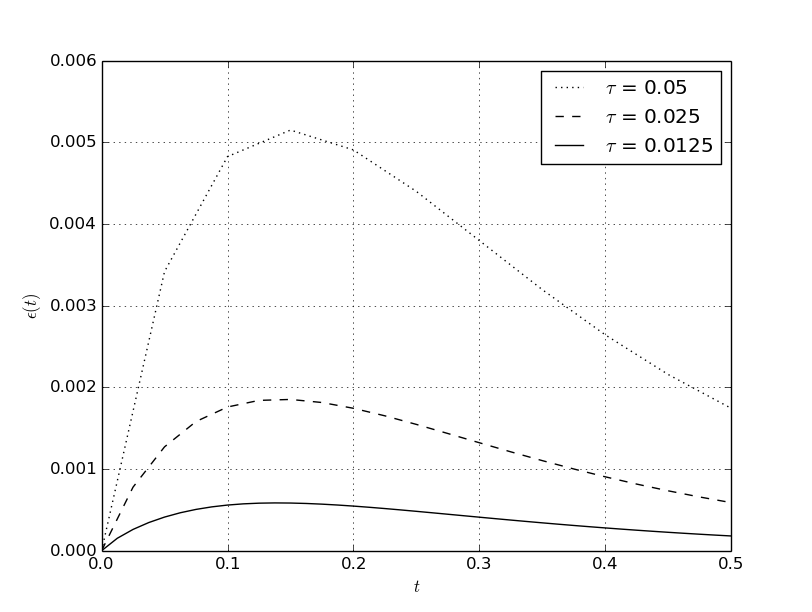}
	\caption{SM-stable factorized scheme: $\sigma = \sqrt{2}$}
	\label{f-9}
  \end{center}
\end{figure}

Similar results obtained using the SM-stable factorized scheme (\ref{25}), (\ref{26})
are presented in Figure~\ref{f-8} and Figure~\ref{f-9}. 
The optimum value ($\sigma = 1/\sqrt{2}$, Figure~\ref{f-8}) provides significantly higher
accuracy in comparison with the case, where we employ the scheme with the twice value of the parameter $\sigma$ 
(Figure~\ref{f-9}).

These results demonstrate the efficiency of  factorized variants
of the SM-stable scheme with the second-order accuracy for solving parabolic boundary value problems.

\clearpage

\end{document}